\documentclass[a4paper,11pt]{article}
\usepackage{latexsym}
\usepackage{amssymb}
\usepackage{enumerate}
\usepackage{amsfonts}
\usepackage{amsmath}
\usepackage[dvips]{graphicx}
\newcommand{\qed}{$\Box$}

\newenvironment{@abssec}[1]{%
    \if@twocolumn

      \section*{#1}%
    \else

      \vspace{.05in}\footnotesize
      \parindent .2in
 {\upshape\bfseries #1. }\ignorespaces
    \fi}

    {\if@twocolumn\else\par\vspace{.1in}\fi}

\newenvironment{keywords}{\begin{@abssec}{\keywordsname}}{\end{@abssec}}

\newenvironment{AMS}{\begin{@abssec}{\AMSname}}{\end{@abssec}}

\newcommand\keywordsname{Key words}
\newcommand\AMSname{AMS subject classifications}
\newcommand\AMname{AMS subject classification}
\newtheorem{theorem}{Theorem}
 \newtheorem{lemma}[theorem]{Lemma}
 \newtheorem{corollary}[theorem]{Corollary}
 
\newtheorem{remark}[theorem]{Remark}
 
\def\qed{\vbox{\hrule height0.6pt\hbox{%
  \vrule height1.3ex width0.6pt\hskip0.8ex
  \vrule width0.6pt}\hrule height0.6pt
 }}

\title{Symmetry problems on stationary isothermic surfaces \\ in Euclidean spaces
\thanks{This research was partially supported by the Grant-in-Aid for Challenging Exploratory Research
($\sharp$ 25610024) of
Japan Society for the Promotion of Science.}}

\author{Shigeru Sakaguchi\thanks{Research Center for Pure and Applied Mathematics,
Graduate School of  Information Sciences, Tohoku
University, Sendai, 980-8579,  Japan.
({\tt sigersak@m.tohoku.ac.jp}).}}
\begin{document}

\maketitle

\begin{abstract}
Let $S$ be a smooth hypersurface properly embedded in $\mathbb R^N$ with $N \geq 3$ and consider its tubular neighborhood $\mathcal N$. We show that, if a heat flow over $\mathcal N$ with appropriate initial and boundary conditions has $S$  as a stationary isothermic surface, then $S$ must have some sort of symmetry.
\end{abstract}


\begin{keywords}
heat equation,  Cauchy problem, initial-boundary value problem, tubular neighborhood, stationary isothermic surface, symmetry.
\end{keywords}

\begin{AMS}
Primary 35K05 ; Secondary  35B40,  35K15, 35K20.
\end{AMS}

\pagestyle{plain}
\thispagestyle{plain}
\markboth{S. SAKAGUCHI}{Symmetry of stationary isothermic surfaces}

\pagestyle{plain}
\thispagestyle{plain}

\section{Introduction}
\label{introduction}

\vskip 2ex
The stationary isothermic surfaces of solutions of the heat equation have been much studied,  and it has been shown that the existence of a stationary isothermic surface forces  the problems  to have some sort of symmetry (see \cite{MPS, MPrStams2006, MSannals2002, MSindiana2007, MSjde2010, MSpoincare2010, MSmmas2013, Sspringer2013}).  A balance law for stationary zeros of temperature  introduced by \cite{MSmathz1999}  plays a key role in the proofs.  To be more precise,  the balance law gives us that  for any pair of points $x$ and $y$ in the stationary isothermic surface the heat contents of two balls centered at $x$ and $y$ respectively with an equal radius are equal for every time.  The above papers always deal with the cases where each ball touches the boundary only at one point eventually. 
Then by studying the initial behavior of the heat content of each ball the authors extract some information of the principal curvatures of the boundary at the touching point. 

We emphasize that in the present paper  we deal with the cases where each ball touches the boundary exactly at two points.
Another new point is to give simply a $C^2$ hypersurface properly embedded in $\mathbb R^N$ as a candidate for a stationary isothermic surface from the beginning.

 Let us establish our setting.  Let $\Omega$ be a $C^2$ domain in $\mathbb R^N$ with $N \geq 3$, whose  boundary $\partial\Omega$ is connected and not necessarily bounded.  Namely,  $\partial\Omega$ is a $C^2$ hypersurface properly embedded in $\mathbb R^N$.  Assume that there exists a number $R > 0$ satisfying:
 \begin{itemize}
 \item[(A-1)]: The principal curvatures $\kappa_1(x), \dots, \kappa_{N-1}(x)$ of $\partial\Omega$  at $x \in \partial\Omega$ with respect to the outward normal direction to $\partial \Omega$ satisfy
 $$
\max_{1 \le j \le N-1} |\kappa_j(x)|  < \frac 1R\ \mbox{ for every } x \in \partial\Omega.
 $$
 \item[(A-2)]: The tubular neighborhood $\mathcal N_R$ of $\partial\Omega$ given by
 $$
 \mathcal N_R = \{ x \in \mathbb R^N : \mbox{ dist}(x, \partial\Omega) < R \},
 $$
  is a  $C^2$ domain in $\mathbb R^N$ and its boundary $\partial  \mathcal N_R$ consists of two connected components $\Gamma_+, \Gamma_-$ each of which is diffeomorphic to $\partial\Omega$.
 \end{itemize}
 Let us introduce two $C^2$ domains $\Omega_+, \Omega_-$  in $\mathbb R^N$ with $\partial\Omega_+ = \Gamma_+,  \partial\Omega_- = \Gamma_-$,  respectively, such that  the three domains $\Omega_+, \Omega_-,  \mathcal N_R$ are disjoint, $\Omega_- \subset \Omega$, and
 $\Omega_+\cup\Omega_-\cup \overline{\mathcal N_R} = \mathbb R^N$.  Denote by ${\mathcal X}_{\Omega_+}, {\mathcal X}_{\Omega_-}$ the characteristic functions of the sets $\Omega_+, \Omega_-$, respectively.
Consider the following initial-boundary value problem for the heat equation:
\begin{eqnarray}
&&u_t =\Delta u\quad\mbox{ in }\ \mathcal N_R\times (0,+\infty), \label{heat equation initial-boundary}
\\
&&u=1 \ \mbox{ on } \partial\mathcal N_R\times (0,+\infty), \label{heat Dirichlet}
\\ 
&&u=0  \ \mbox{ on } \mathcal N_R\times \{0\},\label{heat initial}
\end{eqnarray}
and the Cauchy problem for the heat equation:
\begin{equation}
  u_t =\Delta u\quad\mbox{ in }\  \mathbb R^N\times (0,+\infty) \ \mbox{ and }\ u\ ={\mathcal X}_{\Omega_+}+{\mathcal X}_{\Omega_-}\ \mbox{ on } \mathbb R^N\times
\{0\}.\label{heat Cauchy}
\end{equation}
We have the following theorem.
\begin{theorem} 
\label{th:stationary isothermic1}  Let $N =3$ and let $u$ be the unique bounded solution of either  problem \eqref{heat equation initial-boundary}-\eqref{heat initial} or problem \eqref{heat Cauchy}. Assume that there exists a function $a(t)$ satisfying
\begin{equation}
\label{stationary level 1}
u(x,t) = a(t)\ \mbox{ for every } (x,t) \in \partial\Omega \times (0, +\infty).
\end{equation}
Then, $\partial\Omega$ must be either a plane or a sphere, provided at least one of the following conditions is satisfied:
\begin{itemize}
\item[\rm(a)]\ $\partial\Omega$ has an umbilical point  $p \in \partial\Omega$, that is,  $\kappa_1(p)= \kappa_2(p)$.

\item[\rm(b)]\ There exists a sequence of points $\{ p_j\} \subset \partial\Omega$ with $\lim\limits_{j \to \infty} \kappa_1(p_j) = \lim\limits_{j \to \infty} \kappa_2(p_j) \in \mathbb R$.
\end{itemize}
\end{theorem}

\vskip 2ex
When $\partial\Omega$ is bounded, the Hopf-Poincar\'e theorem \cite[Theorem II, p. 113]{H1989} says that the sum of the indices of all the isolated umbilical points equals the Euler  number $\chi(\partial\Omega) (= 2 -2\times \mbox{ genus })$ of $\partial\Omega$ and hence if the genus of $\partial\Omega$ does not equal $1$ then $\partial\Omega$ must have at least one umbilical point. Therefore we have the following direct corollary. 
\begin{corollary} 
\label{cor:stationary isothermic1} Let $N =3$ and let $u$ be the unique bounded solution of either  problem \eqref{heat equation initial-boundary}-\eqref{heat initial} or problem \eqref{heat Cauchy}. Assume that \eqref{stationary level 1} holds for some function $a(t)$. Then, if  $\partial\Omega$ is bounded and the genus of $\partial\Omega$ does not equal $1$, $\partial\Omega$ must be a sphere.
\end{corollary}

\vskip 2ex
We next consider the following initial-boundary value problem for the heat equation:
\begin{eqnarray}
&&u_t =\Delta u\quad\mbox{ in }\ \mathcal N_R\times (0,+\infty),\label{heat equation plus-minus initial-boundary}
\\
&&u=1 \ \mbox{ on } \Gamma_+\times
(0,+\infty),\label{heat Dirichlet plus}
\\
&&u=-1 \ \mbox{ on } \Gamma_-\times
(0,+\infty),\label{heat Dirichlet-minus}
 \\ 
&&u=0  \ \mbox{ on } \mathcal N_R\times \{0\},\label{heat initial zero}
\end{eqnarray}
and the Cauchy problem for the heat equation:
\begin{equation}
\label{heat Cauchy plus-minus}
 u_t =\Delta u\quad\mbox{ in }\  \mathbb R^N\times (0,+\infty) \ \mbox{ and }\ u\ ={\mathcal X}_{\Omega_+} -{\mathcal X}_{\Omega_-} \ \mbox{ on } \mathbb R^N\times
\{0\}.
\end{equation}
 Then we have
\begin{theorem} 
\label{th:stationary isothermic2}  Let $N \geq 3$ and let  $u$ be the unique bounded solution of either  problem \eqref{heat equation plus-minus initial-boundary}-\eqref{heat initial zero} or problem \eqref{heat Cauchy plus-minus}.
Assume that \eqref{stationary level 1} holds for some function $a(t)$. 
Then:
\begin{itemize}
\item[\rm(1)]\ If $\partial\Omega$ is bounded, 
$\partial\Omega$ must be a sphere. 
\item[\rm(2)]\ If $N=3$ and $\partial\Omega$ is an entire graph over  $\mathbb R^2$,  $\partial\Omega$ must be a plane.
\end{itemize}
\end{theorem}

\vskip 2ex
By using the asymptotic formula of the heat content $\int_{B_R(x)} u(z,t)\ dz$ of an open ball $B_R(x)$ with radius $R > 0$ centered at $ x \in \partial\Omega$ as $t \to +0$ introduced in \cite{MSedinburgh2007} together with the balance law given in \cite{MSmathz1999}, we prove Theorems \ref{th:stationary isothermic1} and \ref{th:stationary isothermic2}.  Moreover Aleksandrov's sphere theorem and Bernstein's theorem for  the minimal surface equation are needed to prove Theorem \ref{th:stationary isothermic2}. In sections \ref{section2} and \ref{section3},  we prove Theorems \ref{th:stationary isothermic1} and \ref{th:stationary isothermic2}, respectively. The final section \ref{section4} gives several remarks and problems.

\setcounter{equation}{0}
\setcounter{theorem}{0}

\section{Proof of Theorem \ref{th:stationary isothermic1}}
\label{section2}

The proofs of Theorems \ref{th:stationary isothermic1} and \ref{th:stationary isothermic2} have common ingredients. Therefore  we begin with general dimensions $N$ for later use,  although Theorem \ref{th:stationary isothermic1} assumes that $N=3$. 

Let $u$ be the unique bounded solution of either  problem \eqref{heat equation initial-boundary}-\eqref{heat initial} or problem \eqref{heat Cauchy}.  Denote by $u^\pm=u^\pm(x,t)$  the unique bounded solutions of the initial-boundary value problems for the heat equation:
\begin{eqnarray}
&&u_t =\Delta u\quad\mbox{ in }\ \left(\mathbb R^N\setminus  \overline{\Omega_\pm}\right)\times (0,+\infty), \label{heat equation initial-boundary pm}
\\
&&u=1 \ \mbox{ on } \Gamma_\pm\times (0,+\infty), \label{heat Dirichlet pm}
\\ 
&&u=0  \ \mbox{ on } \ \left(\mathbb R^N\setminus  \overline{\Omega_\pm}\right)\times \{0\},\label{heat initial pm}
\end{eqnarray}
respectively, 
or of the Cauchy problems for the heat equation:
\begin{equation}
  u_t =\Delta u\quad\mbox{ in }\  \mathbb R^N\times (0,+\infty) \ \mbox{ and }\ u\ ={\mathcal X}_{\Omega_\pm}\ \mbox{ on } \mathbb R^N\times
\{0\},\label{heat Cauchy pm}
\end{equation}
respectively. Notice that $u = u^+ + u^-$  when $u$ is the solution of problem \eqref{heat Cauchy}.
Then, by a result of Varadhan \cite{Vcpam1967}(see also \cite[Theorem A, p. 2024]{MSmmas2013}), we see that
\begin{equation}
\label{varadhan}
-4t\log\left(u^\pm(x,t)\right) \to \mbox{ dist}(x, \Gamma_\pm)^2\ \mbox{ as } t \to +0
\end{equation}
 uniformly on every compact sets in $\mathbb R^N\setminus  \overline{\Omega_\pm}$.

By the assumptions (A-1) and (A-2),  every point $x \in \partial\Omega$ determines two points $x_+ \in \Gamma_+$ and $x_-\in \Gamma_-$ satisfying
$$
\partial B_R(x) \cap \Gamma_+ = \{ x_+\}\ \mbox{ and } \partial B_R(x) \cap \Gamma_- = \{ x_-\},
$$
respectively.  Moreover,  by letting $\kappa^\pm_1(x_\pm),\dots,\kappa^\pm_{N-1}(x_\pm)$ denote the principal curvatures of $\Gamma_\pm$ at $x_\pm$ with respect to the inward normal direction to $\partial \mathcal N_R$,  respectively, we observe that 
\begin{equation}
\label{curvatures relations}
1- R\kappa^+_j(x_+) = \frac 1{1-R\kappa_j(x)} > 0\ \mbox{ and }\ 1- R\kappa^-_j(x_-) =  \frac 1{1+R\kappa_j(x)} >0
\end{equation}
 for every $ x \in \partial\Omega$ and every $j = 1, \dots, N-1$. 
 
 On the other hand, it follows from the balance law (see \cite[Theorem 4, p. 704]{MSmathz1999} or \cite[Theorem 2.1, pp. 934-935]{MSannals2002}) that \eqref{stationary level 1} gives
\begin{equation}
\label{balance law}
\int_{B_R(x)} u(z,t)\ dz = \int_{B_R(y)} u(z,t)\ dz\ \mbox{ for } t >0
\end{equation}
for every $x,y \in \partial\Omega$.  Moreover, by virtue of \eqref{curvatures relations},  an asymptotic formula given by \cite{MSedinburgh2007} (see also
\cite[Theorem B, pp. 2024-2025]{MSmmas2013}) yields that
\begin{equation}
\label{asymptotics and curvatures}
\lim_{t\to +0}t^{-\frac{N+1}4 }\!\!\!\int\limits_{B_R(x)}\! u^\pm(z,t)\ dz=
c(N)\left\{\prod\limits_{j=1}^{N-1}\left[\frac 1R - \kappa^\pm_j(x_\pm)\right]\right\}^{-\frac 12},
\end{equation}
respectively. Here,  $c(N)$ is a positive constant depending only on $N$  and of course $c(N)$ depends on the problems  \eqref{heat equation initial-boundary pm}-\eqref{heat initial pm}  or \eqref{heat Cauchy pm}.  Then we have
\begin{lemma}
\label{le: curvature formula} 
Let $u$ be the unique bounded solution of either  problem \eqref{heat equation initial-boundary}-\eqref{heat initial} or problem \eqref{heat Cauchy}.  Assume that \eqref{stationary level 1} holds for some function $a(t)$. Then there exists a constant $c > 0$ satisfying
\begin{equation}
\label{hyperbolic PDE}
\left\{\prod\limits_{j=1}^{N-1}(1-R\kappa_j(x))\right\}^{\frac 12} + \left\{\prod\limits_{j=1}^{N-1}(1+R\kappa_j(x))\right\}^{\frac 12} = c\ \mbox{ for every } x \in \partial\Omega,
\end{equation}
where $\kappa_1(x),\dots,\kappa_{N-1}(x)$ denote the principal curvatures of $\partial\Omega$ given in {\rm (A-1)}. 
\end{lemma}

\noindent
{\it Proof.\ } Let $u$ be the unique bounded solution of problem \eqref{heat Cauchy}.  Then we have that $u = u^+ + u^-$.
Hence, combining \eqref{balance law} with \eqref{asymptotics and curvatures} yields that there exists a constant $c > 0$ satisfying
\begin{equation}
\label{sum of two points}
\left\{\prod\limits_{j=1}^{N-1}\left(1 - R\kappa^+_j(x_+)\right)\right\}^{-\frac 12} + \left\{\prod\limits_{j=1}^{N-1}\left(1 - R\kappa^-_j(x_-)\right)\right\}^{-\frac 12} = c 
\end{equation}
for every $ x \in \partial\Omega$. Therefore \eqref{curvatures relations} gives the conclusion.

 Let $u$ be the solution of problem \eqref{heat equation initial-boundary}-\eqref{heat initial}. It follows from the comparison principle that
 $$
 \max\{ u^+, u^-\} \le u \le u^+ + u^-\ \ \mbox{ in } \mathcal N_R \times (0,\infty).
 $$
Therefore, in view of  \eqref{varadhan} and \eqref{asymptotics and curvatures}, we notice that for every $x \in \partial\Omega$ 
\begin{eqnarray*}
&&c(N)\left\{\prod\limits_{j=1}^{N-1}\left[\frac 1R - \kappa^+_j(x_+)\right]\right\}^{-\frac 12} + 
c(N)\left\{\prod\limits_{j=1}^{N-1}\left[\frac 1R - \kappa^-_j(x_-)\right]\right\}^{-\frac 12}\\
&& = \lim_{t \to +0} t^{-\frac {N+1}4}\int_{B_R(x)}u^+(z,t)\ dz + \lim_{t \to +0} t^{-\frac {N+1}4}\int_{B_R(x)}u^-(z,t)\ dz\\
&& =\lim_{t \to +0} t^{-\frac {N+1}4}\int_{B_R(x)\setminus\Omega}u^+(z,t)\ dz + \lim_{t \to +0} t^{-\frac {N+1}4}\int_{B_R(x)\cap\Omega}u^-(z,t)\ dz\\
&&= \lim_{t \to +0} t^{-\frac {N+1}4}\int_{B_R(x)\setminus\Omega}u(z,t)\ dz + \lim_{t \to +0} t^{-\frac {N+1}4}\int_{B_R(x)\cap\Omega}u(z,t)\ dz\\
&&=\lim_{t \to +0} t^{-\frac {N+1}4}\int_{B_R(x)}u(z,t)\ dz.
\end{eqnarray*}
Hence, with the aid of \eqref{balance law},  we obtain \eqref{sum of two points} which yields the  conclusion by \eqref{curvatures relations}.  \qed

\vskip 2ex
\noindent
{\bf Proof of Theorem \ref{th:stationary isothermic1}: }
  Set $N=3$ in \eqref{hyperbolic PDE}.  With the aid of the arithmetic-geometric mean inequality, we obtain from \eqref{hyperbolic PDE} that
$$
c = \sqrt{(1-R\kappa_1)(1-R\kappa_2)} + \sqrt{(1+R\kappa_1)(1+R\kappa_2)} \le \frac {2-R(\kappa_1+\kappa_2)}2+ \frac {2+R(\kappa_1+\kappa_2)}2 = 2
$$
where $\kappa_j=\kappa_j(x)$ with $j = 1, 2$. By the assumption, $\partial\Omega$ has an umbilical point  $p \in \partial\Omega$, that is,  $\kappa_1(p)= \kappa_2(p)$, or   there exists a sequence of points $\{ p_j\} \subset \partial\Omega$ with $\lim\limits_{j \to \infty} \kappa_1(p_j) = \lim\limits_{j \to \infty} \kappa_2(p_j) \in \mathbb R$. Then
we conclude that $c = 2$ and the equality holds in the above inequality. Hence $\kappa_1 = \kappa_2$ on $\partial\Omega$, that is, $\partial\Omega$ is called totally umbilical.
Thus from classical results in differential geometry $\partial\Omega$ must be either a plane or a sphere(see \cite[Remark, p. 124]{H1989} or \cite[Theorem 3.30, p. 84]{MR2005gradst69AMS} for instance). \qed

\setcounter{equation}{0}
\setcounter{theorem}{0}

\section{Proof of Theorem \ref{th:stationary isothermic2}}
\label{section3}
 Let us use the auxiliary functions $u^\pm=u^\pm(x,t)$ given in section \ref{section2}. We begin with the following lemma:
\begin{lemma}
\label{le: curvature formula 2} 
Let  $u$ be the unique bounded solution of either  problem \eqref{heat equation plus-minus initial-boundary}-\eqref{heat initial zero} or problem \eqref{heat Cauchy plus-minus}.  Assume that \eqref{stationary level 1} holds for some function $a(t)$. 
 Then there exists a constant $c$ satisfying
\begin{equation}
\label{elliptic PDE}
\left\{\prod\limits_{j=1}^{N-1}(1-R\kappa_j(x))\right\}^{\frac 12} - \left\{\prod\limits_{j=1}^{N-1}(1+R\kappa_j(x))\right\}^{\frac 12} = c\ \mbox{ for every } x \in \partial\Omega,
\end{equation}
where $\kappa_1(x),\dots,\kappa_{N-1}(x)$ denote the principal curvatures of $\partial\Omega$ given in {\rm (A-1)}. 
\end{lemma}

\noindent
{\it Proof.\ } Let  $u$ be the solution of problem \eqref{heat Cauchy plus-minus}. Then we have that $u = u^+ - u^-$. Therefore the conclusion follows from the same argument as in the proof of Lemma \ref{le: curvature formula}.

Let  $u$ be the solution of problem \eqref{heat equation plus-minus initial-boundary}-\eqref{heat initial zero}.
It follows from the comparison principle that
$$
\max\{-u^-, u^+-2u^-\} \le u \le \min\{u^+, 2u^+-u^-\}\ \mbox{ in } \mathcal N_R \times (0,\infty).
$$
With the aid of these inequalities, in view of  \eqref{varadhan} and \eqref{asymptotics and curvatures}, by carrying out calculations similar to those in the proof of Lemma \ref{le: curvature formula} for every $x \in \partial\Omega$, we can reach the conclusion. \qed

\vskip 2ex
\noindent
{\bf Proof of Theorem \ref{th:stationary isothermic2}: }
  Set 
$$
- \Phi(\kappa_1,\dots,\kappa_{N-1}) = \mbox{ the left-hand side of \eqref{elliptic PDE}.}
$$ 
Then we have that $\frac {\partial \Phi}{\partial\kappa_j} > 0$ for $j =1,\dots, N-1$.  Therefore,  by introducing local coordinates, the condition $\Phi(\kappa_1,\dots,\kappa_{N-1}) =$ constant
on the surface $\partial\Omega$ can be converted into a second order partial differential
equation which is of elliptic type. Hence, if $\partial\Omega$ is bounded, then
$\partial\Omega$ must be a sphere by Aleksandrov's sphere theorem \cite{A1958}.  Thus proposition (1) is proved.

Let us proceed to proposition (2).  Set $N=3$ in \eqref{elliptic PDE}. Then
\begin{equation}
\label{Nequals3}
 \sqrt{(1-R\kappa_1)(1-R\kappa_2)} - \sqrt{(1+R\kappa_1)(1+R\kappa_2)} =c,
 \end{equation}
where $\kappa_j=\kappa_j(x)$ with $j = 1, 2$, and hence
\begin{equation}
\label{mean curvature}
-4RH = c\left( \sqrt{(1-R\kappa_1)(1-R\kappa_2)} + \sqrt{(1+R\kappa_1)(1+R\kappa_2)}\right),
\end{equation}
where $H = \frac 12(\kappa_1 +\kappa_2)$ is the mean curvature of $\partial\Omega$. We distinguish three cases: 
$$
{\rm (i)}\  c=0,\ {\rm (ii)}\  c > 0,\ {\rm (iii)}\ c< 0.
$$
In case (i),  by \eqref{mean curvature} we have $H=0$ on $\partial\Omega$ and hence  $\partial\Omega$ is the minimal entire graph of a function over $\mathbb R^2$. Therefore, by Bernstein's theorem for the minimal surface equation,  $\partial\Omega$ must be a plane. This gives the conclusion desired. (See \cite{GT1983, G1984} for Bernstein's theorem.) In case (ii), by \eqref{mean curvature} we have $H<0$ on $\partial\Omega$.  Suppose that there exists a sequence of points $\{ p_n\}$ with $\lim\limits_{n \to \infty}H(p_n) =0$.
Since $R\kappa_1(p_n), R\kappa_2(p_n) \in [-1, 1]$, by the Bolzano-Weierstrass theorem, by taking a subsequence if necessary, we may assume that $\{R\kappa_1(p_n)\},  \{R\kappa_2(p_n)\}$ converge to numbers $\alpha, -\alpha$, respectively,  for some $\alpha \in [-1,1]$.  Hence by \eqref{Nequals3} we get $c =0$ which is a contradiction. Therefore, there exists a number $\delta >0$ such that
$$
H \le -\delta\ \mbox{ on } \partial\Omega,
$$
which contradicts the fact  that $\partial\Omega$ is an entire graph over  $\mathbb R^2$ with the aid of the divergence theorem as in the proof of  \cite[Theorem 3.3,  pp. 2732--2733]{MSindiana2007}. The remaining case (iii) can be dealt with in a similar manner.  Thus proposition (2) is proved. \qed

\begin{remark}
\label{hyperbolic type} In section {\rm \ref{section2}} we did not use the same argument as in section {\rm \ref{section3}},  for by introducing local coordinates, the condition \eqref{hyperbolic PDE} on the surface $\partial\Omega$  can not be converted into a second order partial differential equation which is of elliptic type.
\end{remark}

\setcounter{equation}{0}
\setcounter{theorem}{0}

\section{Concluding Remarks and Problems}
\label{section4}

In this final section, we mention several remarks and problems.  

Concerning Theorem \ref{th:stationary isothermic1}, spherical cylinders satisfy the assumption \eqref{stationary level 1}. Therefore, as in \cite{MPS}, a theorem including a spherical cylinder as a conclusion is expected.
Corollary \ref{cor:stationary isothermic1} excludes closed surfaces with genus 1, but this might be technical. Concerning Theorem \ref{th:stationary isothermic2}, right helicoids satisfy the assumption \eqref{stationary level 1}. Therefore,  a theorem including a right helicoid as a conclusion is expected. 

Let us set $N=3$ both in \eqref{hyperbolic PDE} and in \eqref{elliptic PDE} and assume that $\partial\Omega$ is a minimal surface properly embedded in $\mathbb R^3$.  Then \eqref{hyperbolic PDE} yields that the Gauss curvature is constant and hence $\partial\Omega$ must be a plane. On the other hand, \eqref{elliptic PDE} holds true for every minimal surface by setting $c = 0$. 

Concerning technical points in the theory of partial differential equations, \eqref{hyperbolic PDE} is not of elliptic type but \eqref{elliptic PDE} is of elliptic type, as is mentioned in section \ref{section3}.
Therefore, for \eqref{elliptic PDE} in general dimensions, Liouville-type theorems characterizing hyperplanes are expected as in 
\cite{MSjde2010, Sspringer2013}.

\begin{footnotesize}

\end{footnotesize}
\end{document}